\documentclass[11pt,letterpaper,reqno]{amsart}

\usepackage{amssymb,amsmath,amsthm,amsfonts}
\usepackage{bbm}

\usepackage{bookmark}
\usepackage{hyperref}
\hypersetup{pdfstartview={FitH}}

\newtheorem{theorem}{Theorem}
\newtheorem{problem}[theorem]{Problem}

\numberwithin{equation}{section}

\begin{document}

\title[Convergence of ergodic--martingale paraproducts]{Convergence of ergodic--martingale paraproducts}

\author[Vjekoslav Kova\v{c}]{Vjekoslav Kova\v{c}}
\address{Vjekoslav Kova\v{c}, Department of Mathematics, Faculty of Science, University of Zagreb, Bijeni\v{c}ka cesta 30, 10000 Zagreb, Croatia}
\email{vjekovac@math.hr}

\author[Mario Stip\v{c}i\'{c}]{Mario Stip\v{c}i\'{c}}
\address{Mario Stip\v{c}i\'{c}, Department of Mathematics, Faculty of Science, University of Zagreb, Bijeni\v{c}ka cesta 30, 10000 Zagreb, Croatia}
\email{mstipcic@math.hr}




\begin{abstract}
In this note we introduce a sequence of bilinear operators that unify ergodic averages and backward martingales in a nontrivial way. We establish its convergence in a range of $\textup{L}^p$-norms and leave its a.s.\@ convergence as an open problem. This problem shares some similarities with a well-known unresolved conjecture on a.s.\@ convergence of double ergodic averages with respect to two commuting transformations.
\end{abstract}

\maketitle


\section{Introduction}
There are many similarities in the behaviors of ergodic averages and (forward or backward) martingales. Back in 1950 they inspired Kakutani \cite{K52:open} to formulate an open-ended problem of finding a single concept that generalizes both of these notions. He was primarily looking for ``a general theorem which contains both the maximal ergodic theorem and the martingale theorem'' (a quote from \cite{K52:open}), and both of these are results on convergence almost surely. However, one can understand his question in a broader sense, by also considering other modes of convergence. Kakutani's question was answered in versatile ways by many different authors over the course of the last 70 years. The most notable unifying theories were developed by Jerison \cite{J59:open}, Rota \cite{R61:open}, A. and C. Ionescu Tulcea \cite{IT63:open}, Petz \cite{P84:open}, Kachurovskii \cite{K98:open}, and Kachurovskii and Vershik \cite{VK99:open}; see the survey by Kachurovskii \cite{K07:open}. It is also interesting to mention a largely forgotten paper of Neveu \cite{N64:open}, who deduced a.s.\@ convergence of backward martingales from the pointwise ergodic theorem for contractions. The question of unifying ergodic averages and martingales still attracts some attention of the mathematical community; see the more recent papers by Podvigin \cite{P10:open,P14:open}, Ganiev and Shahidi \cite{SG12:open}, and Shahidi \cite{S20:open}.

One incentive for writing this note was our wish to approach the aforementioned question of Kakutani via bilinear operators and in the spirit of classical harmonic analysis. To our surprise, already the simplest nontrivial bilinear objects formed by ergodic averages and discrete martingales turned out to be somewhat involved.

Let $(\Omega,\mathcal{F},\mathbb{P})$ be a probability space. The conditional expectation operator $\mathbb{E}(\cdot|\mathcal{G})$ is associated with each $\sigma$-algebra $\mathcal{G}\subseteq\mathcal{F}$.
A \emph{backward filtration} of $(\Omega,\mathcal{F})$ is any decreasing sequence $(\mathcal{G}_n)_{n=0}^{\infty}$ of $\sigma$-algebras such that $\mathcal{G}_0=\mathcal{F}$. A \emph{backward martingale} with respect to that filtration is any sequence $(f_n)_{n=0}^{\infty}$ of real-valued functions in $\textup{L}^1(\Omega,\mathcal{F},\mathbb{P})$ such that each $f_n$ is $\mathcal{G}_n$-measurable and $\mathbb{E}(f_{n}|\mathcal{G}_{n+1})=f_{n+1}$ a.s.\@ for each index $n$.
Any given function $g\in\textup{L}^p(\Omega,\mathcal{F},\mathbb{P})$, $1\leq p<\infty$, gives rise to one such backward martingale, namely $(E_n g)_{n=0}^{\infty}$ defined by
\begin{equation}\label{eq:defmartingale}
E_n g :=\mathbb{E}(g|\mathcal{G}_n)
\end{equation}
for every nonnegative integer $n$. As a byproduct of the proof of Doob's martingale convergence theorem \cite{D40:mart} this sequence converges in the $\textup{L}^p$-norm and a.s.\@ as $n\to\infty$ and the limit can be identified as $\mathbb{E}(g|\cap_{n=0}^{\infty}\mathcal{G}_n)$. Conversely, all backward martingales can be obtained via the above construction. Proofs of these results can be found in many graduate level textbooks on the probability theory basics, such as the one by Durrett \cite{D10:prob}.

Suppose that a transformation $T\colon\Omega\to\Omega$ is $(\mathcal{F},\mathcal{F})$-\emph{measurable}, i.e.\@ $T^{-1}E\in\mathcal{F}$ for every $E\in\mathcal{F}$, and \emph{measure-preserving}, i.e.\@ $\mathbb{P}(T^{-1}E)=\mathbb{P}(E)$ for every $E\in\mathcal{F}$. For a positive integer $k$ we denote by $T^k$ the $k$-th iterate of $T$, i.e., the $k$-fold composition $T\circ\cdots\circ T$, while $T^0$ is interpreted as the identity on $\Omega$. Central objects in classical ergodic theory are the \emph{(Ces\`{a}ro) ergodic averages} $A_N f$ of an $\mathcal{F}$-measurable function $f\colon\Omega\to\mathbb{R}$ with respect to the iterates of $T$. Namely,
\begin{equation}\label{eq:defergaver}
A_N f := \frac{1}{N}\sum_{k=0}^{N-1} f\circ T^k
\end{equation}
for every positive integer $N$.
If $f\in\textup{L}^p(\Omega,\mathcal{F},\mathbb{P})$, $1\leq p<\infty$, then von Neumann's mean ergodic theorem \cite{vN32:erg} establishes convergence of \eqref{eq:defergaver} as $N\to\infty$ in the $\textup{L}^p$-norm, while Birkhoff's pointwise ergodic theorem \cite{B31:pt} gives its convergence a.s. This time the limit can be identified as the conditional expectation of $f$ with respect to the so-called invariant $\sigma$-algebra.
Over the years it became clear that it is actually lacunary subsequences of $(A_N f)_{N=1}^\infty$ that emulate properties of backward martingales, rather than the whole sequence itself. However, it is an easy exercise to deduce convergence of $(A_N f)_{N=1}^\infty$, either in $\textup{L}^p$ or a.s., solely from the fact that $(A_{\lfloor a^n\rfloor}f)_{n=0}^{\infty}$ converges for each $a\in(1,\infty)$; see Appendix of the paper by Frantzikinakis, Lesigne, and Wierdl \cite{FLW16:random}. Here $\lfloor x\rfloor$ denotes the greatest integer not exceeding a given real number $x$.

In what follows we will need some compatibility between the backward filtration $(\mathcal{G}_n)_{n=0}^{\infty}$ and the transformation $T$. We find it natural to impose the \emph{commutativity condition}, i.e., we require that the operators $f\mapsto f\circ T$ and $f\mapsto \mathbb{E}(f|\mathcal{G}_n)$ commute. This means
\begin{equation}\label{eq:commcond}
\mathbb{E}(f\circ T|\mathcal{G}_n) = \mathbb{E}(f|\mathcal{G}_n) \circ T
\end{equation}
for each nonnegative integer $n$ and each function $f\in\textup{L}^1(\Omega,\mathcal{F},\mathbb{P})$. Condition \eqref{eq:commcond} is equivalent to asking that $E_n A_N = A_N E_n$ holds for every $n\geq 0$ and $N\geq 1$, as an equality of operators on $\textup{L}^1(\Omega,\mathcal{F},\mathbb{P})$. A typical example when \eqref{eq:commcond} is satisfied is the case of a bijective measure-preserving transformation $T$ such that both $T$ and $T^{-1}$ are $(\mathcal{G}_n,\mathcal{G}_n)$-measurable for each $n$; we leave this fact as an exercise to the reader. The commutativity condition \eqref{eq:commcond} has already been utilized by Podvigin \cite{P10:open,P14:open}, even though in a different context than ours.

Throughout the paper we will assume that a number $a\in(1,\infty)$ is fixed. We will also be working with triples of exponents $p,q,r\in[1,\infty)$ satisfying the so-called \emph{H\"{o}lder scaling}, i.e., $1/r = 1/p + 1/q$.
A completely trivial construction combining the subsequence of ergodic averages $(A_{\lfloor a^n\rfloor} f)_{n=0}^{\infty}$ and the backward martingale $(E_n g)_{n=0}^{\infty}$ is their pointwise product
\begin{equation}\label{eq:prod}
(A_{\lfloor a^n\rfloor} f) (E_n g).
\end{equation}
Take exponents $p,q,r$ as above and suppose that $f\in\textup{L}^p(\Omega,\mathcal{F},\mathbb{P})$ and $g\in\textup{L}^q(\Omega,\mathcal{F},\mathbb{P})$. Simply by H\"{o}lder's inequality and the aforementioned classical results we see that the sequence of products \eqref{eq:prod} converges in the $\textup{L}^r$-norm and a.s.\@ as $n\to\infty$. By taking either $f$ or $g$ to be constantly equal to $1$ we recover convergence of the backward martingale \eqref{eq:defmartingale} or convergence of the ergodic averages \eqref{eq:defergaver} alone. Needless to say, this generalization of the two concepts is not sufficiently challenging.

The \emph{ergodic--martingale paraproduct} (with respect to $T$ and $(\mathcal{G}_n)_{n=0}^{\infty}$) is the sequence\linebreak 
$(\Pi_n^{\textup{em}})_{n=1}^{\infty}$ of bilinear operators $(f,g)\mapsto\Pi_n^{\textup{em}}(f,g)$ given by
\begin{equation}\label{eq:pprodem}
\Pi_n^{\textup{em}}(f,g) := \sum_{k=0}^{n-1} (A_{\lfloor a^k\rfloor} f) (E_{k+1} g - E_{k} g)
\end{equation}
for every positive integer $n$ and functions $f\in\textup{L}^p(\Omega,\mathcal{F},\mathbb{P})$, $g\in\textup{L}^q(\Omega,\mathcal{F},\mathbb{P})$.
One can equally well define a complementary object, the \emph{martingale--ergodic paraproduct}, which is the sequence of bilinear operators $(\Pi_n^{\textup{me}})_{n=1}^{\infty}$ defined as
\begin{equation}\label{eq:pprodme}
\Pi_n^{\textup{me}}(f,g) := \sum_{k=0}^{n-1} (A_{\lfloor a^{k+1}\rfloor} f - A_{\lfloor a^{k}\rfloor} f) (E_{k+1} g).
\end{equation}
Summation by parts gives
\begin{equation}\label{eq:pprodsum}
\Pi_n^{\textup{em}}(f,g) + \Pi_n^{\textup{me}}(f,g) = (A_{\lfloor a^n\rfloor} f) (E_n g) - f g.
\end{equation}
Therefore the sum of the two paraproducts is a trivial object, known to converge both in the $\textup{L}^r$-norm and a.s.
In contrast with that, it is not clear if either \eqref{eq:pprodem} or \eqref{eq:pprodme} by itself converges in any sense, for any choice of the exponents $p$ and $q$.

Let us shortly justify the term ``paraproduct'' in relation with \eqref{eq:pprodem} and \eqref{eq:pprodme}.
An expository note by B\'{e}nyi, Maldonado, and Naibo \cite{BMN10:pprod} lists several examples of bilinear objects in the harmonic analysis literature that all deserve to be named ``paraproducts.'' In probability theory, paraproducts of two martingales appear as variants of Burkholder's martingale transforms \cite{Bur66:marttran}. Development of their theory, parallel to the one of analytical paraproducts, was initiated by Ba\~{n}uelos and Bennett \cite{BB88:martpprod} (in continuous time) and Chao and Long \cite{CL92:martpprod} (in discrete time); also see more recent papers \cite{KS17:mart,KS15:mart,KZK19:mart}. Here we take even more liberty with usage of the word, as we only have two alike objects that add up to the pointwise product \eqref{eq:prod}, modulo the trivial term $fg$.

Now we formulate the main result of this note.

\begin{theorem}\label{thm:normconvergence}
Take $a\in(1,\infty)$ and suppose that $p,q\in[4/3,4]$, $r\in[1,4/3]$ satisfy the H\"{o}lder scaling. For any functions $f\in\textup{L}^p(\Omega,\mathcal{F},\mathbb{P})$ and $g\in\textup{L}^q(\Omega,\mathcal{F},\mathbb{P})$ the sequences $(\Pi_n^{\textup{em}}(f,g))_{n=1}^{\infty}$ and $(\Pi_n^{\textup{me}}(f,g))_{n=1}^{\infty}$ given by \eqref{eq:pprodem} and \eqref{eq:pprodme}, respectively, converge in the $\textup{L}^r$-norm.
\end{theorem}

By monotonicity of the $\textup{L}^r$-norms on a probability space, one can freely lower the exponent $r$, keeping $p$ and $q$ fixed. However, we prefer to stay within the H\"{o}lder scaling, which is present in the analysis on $\mathbb{R}$ or $\mathbb{R}^d$, where the underlying measure space is only $\sigma$-finite. Moreover, we have formulated the theorem for $r\geq 1$ only, even though some estimates mapping below $\textup{L}^1$ will appear in the proof as useful intermediate steps of multilinear interpolation. Even with these constraints it is very unlikely that the range of exponents $p,q,r$ in Theorem~\ref{thm:normconvergence} is the largest possible one.

For the proof of Theorem~\ref{thm:normconvergence} we will need to invoke a combination of quite a few results from the literature, many of which became available only recently. Therefore, this note does not provide any brand new tricks or techniques, but rather assembles the existing ingredients into a complete proof of the above result.

Properties of \eqref{eq:pprodem} and \eqref{eq:pprodme} still seem to be far from completely understood. Proving a result on their a.s.\@ convergence would certainly be more satisfactory and more in line with Kakutani's question. However, we do not find any techniques in the literature that could solve this problem and we leave it as an open question, hoping that it would attract attention of mathematicians with diverse backgrounds.

\begin{problem}\label{prob:asconvergence}
Prove or disprove that for every pair of functions $f,g\in\textup{L}^{\infty}(\Omega,\mathcal{F},\mathbb{P})$ and every $a\in(1,\infty)$ the sequences $(\Pi_n^{\textup{em}}(f,g))_{n=1}^{\infty}$ and $(\Pi_n^{\textup{me}}(f,g))_{n=1}^{\infty}$ converge a.s.
\end{problem}

Problem~\ref{prob:asconvergence} has a similar flavor as a longstanding conjecture in ergodic theory that we are about to formulate.
Now let $S,T\colon\Omega\to\Omega$ be two arbitrary $(\mathcal{F},\mathcal{F})$-measurable measure-preserving transformations and suppose that they commute, i.e., $ST=TS$. For any two $\mathcal{F}$-measurable functions $f,g\colon\Omega\to\mathbb{R}$ and a positive integer $N$ one can define \emph{double (or bilinear) ergodic averages} $B_N(f,g)$ by
\begin{equation}\label{eq:doubleaver}
B_N(f,g) := \frac{1}{N}\sum_{k=0}^{N-1} (f\circ S^k) (g\circ T^k).
\end{equation}

\begin{problem}\label{prob:openergodic}
Prove or disprove that for every pair of functions $f,g\in\textup{L}^{\infty}(\Omega,\mathcal{F},\mathbb{P})$ the sequence $(B_N(f,g))_{N=1}^{\infty}$ converges a.s.
\end{problem}

Most authors attribute Problem~\ref{prob:openergodic} to either Calder\'{o}n or Furstenberg and it appears on many lists of open problems in ergodic theory, such as the one by Frantzikinakis \cite{Fra16:problems}. In any case, Furstenberg and Katznelson \cite{FK78:msz} were the first to prove any results on multiple ergodic averages associated with mutually commuting transformations. However, they were motivated by applications to additive combinatorics and they did not study convergence.

It is an easy observation (see the appendix of \cite{FLW16:random}) that Problem~\ref{prob:openergodic} would be solved affirmatively if one could only prove that the limit of $B_{\lfloor a^n\rfloor}(f,g)$ as $n\to\infty$ exists a.s.\@ for every $a\in(1,\infty)$ and every pair of $\textup{L}^\infty$ functions $f,g$. On the other hand, $\textup{L}^2$ convergence of the sequence $(B_N(f,g))_{n=1}^{\infty}$ was confirmed by Conze and Lesigne \cite{CL84:L2}. Durcik, \v{S}kreb, Thiele, and one of the present authors \cite{DKST17:nvea} gave an alternative and more quantitative proof of the same fact by showing that, for any given $\varepsilon>0$, this sequence makes $O(\varepsilon^{-2})$ jumps in the $\textup{L}^2$-norm. The proof from \cite{DKST17:nvea} reduces estimates for \eqref{eq:doubleaver} to bounds for non-typical paraproduct-type operators somewhat similar to \eqref{eq:pprodem} and \eqref{eq:pprodme}. Because of this we think that any solution to Problem~\ref{prob:asconvergence} is likely to make progress on Problem~\ref{prob:openergodic} as well. This belief serves as another source of our motivation for writing this note.

\section{Proof of Theorem~\ref{thm:normconvergence}}
We will prove $\textup{L}^r$-convergence of \eqref{eq:pprodem} only, while the convergence of \eqref{eq:pprodme} will then be an immediate consequence of \eqref{eq:pprodsum}.
Our main task is to prove the estimate
\begin{equation}\label{eq:mainestimate}
\Big\| \sum_{k=0}^{n-1} (A_{\lfloor a^k\rfloor} f) (E_{k+1} g - E_{k} g) \Big\|_{\textup{L}^r}
\leq C_{a,p,q,r} \|f\|_{\textup{L}^p} \|g\|_{\textup{L}^q},
\end{equation}
making sure that the constant $C_{a,p,q,r}$ does not depend on $(\mathcal{G}_n)_{n=0}^{\infty}$, $T$, $n$, $f$, or $g$.
Once \eqref{eq:mainestimate} is established, we can take an integer $m$ such that $1\leq m<n$ and apply Estimate \eqref{eq:mainestimate} to the function $E_n g - E_m g$ in place of $g$, which gives
\[ \big\| \Pi_n^{\textup{em}}(f,g) - \Pi_m^{\textup{em}}(f,g) \big\|_{\textup{L}^r}
\leq C_{a,p,q,r} \|f\|_{\textup{L}^p} \|E_n g - E_m g\|_{\textup{L}^q}. \]
The classical results on backward martingales mentioned in the previous section ensure that $(E_n g)_{n=0}^{\infty}$ is a Cauchy sequence in $\textup{L}^{q}(\Omega,\mathcal{F},\mathbb{P})$. This will imply that $(\Pi_n^{\textup{em}}(f,g))_{n=1}^{\infty}$ is a Cauchy sequence in $\textup{L}^{r}(\Omega,\mathcal{F},\mathbb{P})$ and Theorem~\ref{thm:normconvergence} will be established.

As the first reduction in the proof of \eqref{eq:mainestimate}, we can assume that $f$ and $g$ are simple functions, i.e., finite linear combinations of indicator functions of sets from $\mathcal{F}$. The general case then follows by density of the set of simple functions in $\textup{L}^{p}(\Omega,\mathcal{F},\mathbb{P})$ and $\textup{L}^{q}(\Omega,\mathcal{F},\mathbb{P})$, combined with H\"{o}lder's inequality and a trivial fact that the operators $A_N$ and $E_n$ are bounded on the Lebesgue spaces for every $N$ and $n$. This qualitative assumption will be convenient later in the proof.

Second, we can reduce the case of a general $a\in(1,\infty)$ to the particular case $a=2$. Indeed, for a nonnegative integer $l$ let $K(l)$ be the smallest nonnegative integer $k$ such that $\lfloor a^k\rfloor\geq 2^l$. We want to compare the left hand side of \eqref{eq:mainestimate} to
\begin{equation}\label{eq:reduced}
\bigg\| \sum_{l=0}^{\lfloor (n-1)\log_2 a\rfloor} (A_{2^l} f) (E_{K(l+1)} g - E_{K(l)} g) \bigg\|_{\textup{L}^r} ,
\end{equation}
by comparing $A_{\lfloor a^k\rfloor} f$ to $A_{2^l} f$ for indices $k$ such that $K(l)\leq k<K(l+1)$.
For that purpose we need a square-function type estimate for ergodic averages,
\begin{equation}\label{eq:ergodicsquare}
\sup_{N_1<N_2<N_3<\cdots} \Big\| \Big( \sum_{i=1}^{\infty} |A_{N_{i+1}}f - A_{N_{i}}f|^2 \Big)^{1/2} \Big\|_{\textup{L}^p}
\leq C_p \|f\|_{\textup{L}^p},
\end{equation}
where the supremum is taken over all strictly increasing sequences of positive integers $(N_i)_{i=1}^{\infty}$. Bound \eqref{eq:ergodicsquare} for $p\in(1,2]$ is a result by Jones, Ostrovskii, and Rosenblatt \cite[Theorem~2.6]{JOR96:L2}, while the cases $p\in(2,\infty)$ were covered by Jones, Kaufman, Rosenblatt, and Wierdl \cite[Theorem~4.6 for $A_1$]{JKRW98:Lp}.
On the other hand, the bound
\begin{equation}\label{eq:martsquare}
\Big\| \Big( \sum_{k=0}^{\infty} |E_{k+1}g - E_{k}g|^2 \Big)^{1/2} \Big\|_{\textup{L}^q}
\leq C_q \|g\|_{\textup{L}^q}
\end{equation}
for $q\in(1,\infty)$ is simply Burkholder's estimate for the martingale square function \cite{Bur66:marttran}.
Estimating the number of integers $\lfloor a^k\rfloor$ that fall into $[2^l,2^{l+1})$ by at most $\lfloor \log_a 2 \rfloor + 1$, using the Cauchy--Schwarz inequality for discrete sums, and applying \eqref{eq:ergodicsquare} and \eqref{eq:martsquare}, we easily see that the left hand side of \eqref{eq:mainestimate} differs from \eqref{eq:reduced} by at most a constant depending on $a,p,q$ times $\|f\|_{\textup{L}^p}\|g\|_{\textup{L}^q}$.
The sum in \eqref{eq:reduced} is of the same form as the one appearing in \eqref{eq:mainestimate}, except that the backward martingale $(E_i g)_{i=0}^{\infty}$ is sampled at times $K(l)$, which amounts to considering a subsequence of the initial backward filtration $(\mathcal{G}_i)_{i=0}^{\infty}$.

In the third reduction step, we transfer Estimate \eqref{eq:mainestimate}, now specialized to $a=2$, first to the product space $\mathbb{Z}\times\Omega$, and then to $\mathbb{R}\times\Omega$. This procedure will be a variant of the so-called Calder\'{o}n's transference trick \cite{Cal68:transf}, which reduces ergodic averages to averages of real-variable functions.
Suppose that we can prove the estimate
\begin{multline}\label{eq:realest}
\bigg\| \sum_{k=0}^{n-1} \Big( \frac{1}{2^k} \int_{0}^{2^k} F(x+y,\omega) \,\textup{d}y \Big) \Big( \mathbb{E}(G(x,\omega)|\mathcal{G}_{k+1}) - \mathbb{E}(G(x,\omega)|\mathcal{G}_{k}) \Big) \bigg\|_{\textup{L}^r_{(x,\omega)}(\mathbb{R}\times\Omega)} \\
\leq C_{p,q,r} \|F\|_{\textup{L}^p(\mathbb{R}\times\Omega)} \|G\|_{\textup{L}^q(\mathbb{R}\times\Omega)}
\end{multline}
for every pair of simple jointly measurable functions $F,G\colon\mathbb{R}\times\Omega\to\mathbb{R}$.
We introduce $\widetilde{F}\colon\mathbb{Z}\times\Omega\to\mathbb{R}$ as a function along the forward orbits of $T$, defined as $\widetilde{F}(m,\omega) := f(T^m\omega)$ for $m=0,1,\ldots,2^{n+1}-1$ and $0$ otherwise.
Then we take $F(x,\omega)=\sum_{m\in\mathbb{Z}}\widetilde{F}(m,\omega)\mathbbm{1}_{[m,m+1)}(x)$, so that measure-invariance of $T$ gives
\[ \|f\|_{\textup{L}^p}^p = \frac{1}{2^{n+1}} \sum_{m=0}^{2^{n+1}-1} \|f\circ T^m\|_{\textup{L}^p}^p
= \frac{1}{2^{n+1}} \big\|\widetilde{F}\big\|_{\textup{L}^p(\mathbb{Z}\times\Omega)}^p
= \frac{1}{2^{n+1}} \|F\|_{\textup{L}^p(\mathbb{R}\times\Omega)}^p. \]
The same thing can be done with the function $g$, introducing the analogous $\widetilde{G}$ and taking analogously defined $G$.
Measure-invariance of $T$ and the commutativity condition \eqref{eq:commcond} allow us to write
\begin{multline}\label{eq:integerest}
\Big\| \sum_{k=0}^{n-1} \Big( \frac{1}{2^k} \sum_{j=0}^{2^k-1} f\circ T^j \Big) (E_{k+1} g - E_{k} g) \Big\|_{\textup{L}^r}^r \\
\leq \frac{1}{2^n} \bigg\| \sum_{k=0}^{n-1} \Big( \frac{1}{2^k} \sum_{j=0}^{2^k-1} \widetilde{F}(m+j,\omega) \Big) \big( \mathbb{E}(\widetilde{G}(m,\omega)|\mathcal{G}_{k+1}) - \mathbb{E}(\widetilde{G}(m,\omega)|\mathcal{G}_{k}) \big) \bigg\|_{\textup{L}^r_{(m,\omega)}(\mathbb{Z}\times\Omega)}^r.
\end{multline}
It is easy to see that the left hand side of \eqref{eq:realest} differs from the last $\textup{L}^r$-norm in \eqref{eq:integerest} by at most $8\|\widetilde{F}\|_{\textup{L}^p(\mathbb{Z}\times\Omega)}\|\widetilde{G}\|_{\textup{L}^q(\mathbb{Z}\times\Omega)}$.
Therefore, \eqref{eq:mainestimate} will be a consequence of \eqref{eq:realest}, once we prove the latter estimate.
Transference procedure very similar to the one above was performed in \cite[Section~5]{DKST17:nvea}.

The fourth step of the proof reduces averages over intervals in $\mathbb{R}$ to dyadic martingales on $[0,1)$. Before anything else, we observe that \eqref{eq:realest} can be safely rescaled from $[0,2^{n+1})$ to the unit interval by applying $x\mapsto 2^{-n-1}x$ in the first coordinates of $F$ and $G$. This is enabled by the fact that we are working in the H\"{o}lder range of exponents. Jones, Kaufman, Rosenblatt, Wierdl \cite[Theorem~C]{JKRW98:Lp} showed the estimate
\[ \bigg\| \bigg( \sum_{j=0}^{\infty} \Big| 2^{j} \int_{0}^{2^{-j}} h(x+y) \,\textup{d}y - \mathbb{E}(h|\mathcal{D}_j)(x) \Big|^2 \bigg)^{1/2} \bigg\|_{\textup{L}^p_x(\mathbb{R})} \leq C_p \|h\|_{\textup{L}^p(\mathbb{R})} \]
for $p\in(1,\infty)$. Here $(\mathcal{D}_j)_{j=0}^{\infty}$ denotes the standard dyadic forward filtration of $[0,1)$, i.e., $\mathcal{D}_j$ is generated by the intervals $[l2^{-j},(l+1)2^{-j})$; $l=0,1,\ldots,2^{j}-1$. Applying the last estimate in the first variable, invoking \eqref{eq:martsquare} in the second variable, and using the Cauchy--Schwarz inequality for sums in $k$, we easily reduce \eqref{eq:realest} to
\begin{multline*}
\Big\| \sum_{k=0}^{n-1} \mathbb{E}_1(F|\mathcal{D}_{n+1-k}) \big( \mathbb{E}_2(G|\mathcal{G}_{k+1}) - \mathbb{E}_2(G|\mathcal{G}_{k}) \big) \Big\|_{\textup{L}^r([0,1)\times\Omega)} \\
\leq C_{p,q,r} \|F\|_{\textup{L}^p([0,1)\times\Omega)} \|G\|_{\textup{L}^q([0,1)\times\Omega)} .
\end{multline*}
Subscripts in $\mathbb{E}_1$ and $\mathbb{E}_2$ denote in which of the two variables the conditional expectation is taken.
By merely reversing the order of summation in $k$ we can rewrite the last estimate as
\begin{multline}\label{eq:auxmartingaleest}
\Big\| \sum_{k=0}^{n-1} \mathbb{E}_1(F|\mathcal{U}_{k}) \big( \mathbb{E}_2(G|\mathcal{V}_{k+1}) - \mathbb{E}_2(G|\mathcal{V}_{k}) \big) \Big\|_{\textup{L}^r(\Omega_1\times\Omega_2)} \\
\leq C_{p,q,r} \|F\|_{\textup{L}^p(\Omega_1\times\Omega_2)} \|G\|_{\textup{L}^q(\Omega_1\times\Omega_2)} ,
\end{multline}
where $(\mathcal{U}_i)_{i=0}^{\infty}$ and $(\mathcal{V}_i)_{i=0}^{\infty}$ are now forward filtrations of two probability spaces $\Omega_1$ and $\Omega_2$, respectively, i.e., increasing sequences of $\sigma$-algebras on those spaces.
The expression in \eqref{eq:auxmartingaleest} could be called the \emph{martingale--martingale paraproduct}. We emphasize that it is different from the more classical martingale paraproducts appearing in \cite{BB88:martpprod,Bur66:marttran,CL92:martpprod,KS17:mart,KZK19:mart}, because conditional expectations with respect to different filtrations are applied to the functions $F$ and $G$. However, \v{S}kreb and one of the present authors \cite{KS15:mart} have already studied this object to some extent, which is what we find convenient below.

The fifth reduction step is merely a simple observation that we only need to prove \eqref{eq:auxmartingaleest} when all $\sigma$-algebras $\mathcal{U}_i,\mathcal{V}_i$ are finite and, consequently, atomized. Indeed, we have only applied \eqref{eq:auxmartingaleest} with $(\mathcal{U}_i)_{i=0}^{\infty}$ being the shifted dyadic-filtration of $[0,1)$.
To see that each $\mathcal{V}_i$ can be replaced with a countably generated $\sigma$-algebra, we first consider all sets $A_1,\ldots,A_m$ appearing in representations of $G(x,\cdot)$ as simple functions. Indeed, since $g$ was assumed to be simple, there are only finitely many different fibers $G(x,\cdot)$ as $x$ ranges over $[0,1)$. Next, we define $\widetilde{\mathcal{V}}_i$ to be the $\sigma$-algebra generated by the sets $\{\mathbb{E}(\mathbbm{1}_{A_j}|\mathcal{V}_l)>\alpha\}$ for $\alpha\in\mathbb{Q}$, $l\in\{0,1,\ldots,i\}$, $j\in\{1,\ldots,m\}$, and observe that $\mathbb{E}(\mathbbm{1}_{A_j}|\widetilde{\mathcal{V}}_i)=\mathbb{E}(\mathbbm{1}_{A_j}|\mathcal{V}_i)$ a.s.\@ for every indices $i$ and $j$.
Passage from finitely generated to countably generated $\sigma$-algebras is easily performed using convergence theorems for forward martingales; see \cite[Theorem~5.5.7]{D10:prob}.
Moreover, we can even reduce \eqref{eq:auxmartingaleest} to the case when the forward filtrations $(\mathcal{U}_i)_{i=0}^{\infty}$ and $(\mathcal{V}_i)_{i=0}^{\infty}$ progress by always splitting each atom into at most two subatoms. This is easily achieved by inserting intermediate $\sigma$-algebras between $\mathcal{V}_k$ and $\mathcal{V}_{k+1}$.

We are finally done with all of the reductions and turn to providing arguments that establish Estimate \eqref{eq:auxmartingaleest}.
The particular case $p=4$, $q=2$ already exists in the literature. The simplest argument is to invoke \cite[Theorem~1.(a)]{KS15:mart}, which proves exactly that estimate, even for two completely general forward filtrations acting in different coordinates of the product space $\Omega_1\times\Omega_2$. A simpler proof concerning two dyadic filtrations was previously given by one of the present authors \cite[Sections~3\&8]{K12:tp}. The same arguments also apply here after $(\mathcal{U}_i)_{i=0}^{\infty}$ and $(\mathcal{V}_i)_{i=0}^{\infty}$ have been reduced as in the previous paragraph; one only needs to redefine the Haar functions to accommodate for possibly uneven splitting of atoms. This particular case of \eqref{eq:auxmartingaleest} is the most difficult ingredient in our proof; we encourage the reader to go over the details in either \cite{KS15:mart} or \cite{K12:tp}.

At the very end, we want to obtain Estimate \eqref{eq:auxmartingaleest} for the full claimed range of exponents.
One can first use summation by parts to interchange the roles of $(\mathcal{U}_i)_{i=0}^{\infty}$ and $(\mathcal{V}_i)_{i=0}^{\infty}$, and, respectively, also the roles of $F$ and $G$. That way the argument from the previous paragraph also gives Estimate \eqref{eq:auxmartingaleest} for $p=2$, $q=4$, and then multilinear interpolation of $\textup{L}^p$ spaces \cite{GLLZ:interp} yields \eqref{eq:auxmartingaleest} for pairs $(1/p,1/q)$ on the segment joining the points $(1/4,1/2)$ and $(1/2,1/4)$. Afterwards, one uses Gundy's decomposition \cite{G68:decomp} of the forward martingale $(\mathbb{E}(G(\omega,\cdot)|\mathcal{V}_i))_{i=0}^{\infty}$ for each fixed $\omega\in\Omega_1$, in the same way the Calder\'{o}n--Zygmund decomposition was used in \cite[Section~5]{K12:tp}. Our assumptions on the functions and the filtrations come in very handy here, because they resolve any measurability issues. This decomposition gives weak-type bounds at the endpoint $q=1$. Then one uses multilinear interpolation again to recover the part of the claimed range with $p\in[2,4]$. The rest of the desired range is covered using the summation by parts again and repeating the same procedure, this time with $F$ and decomposing in the first coordinate.
This completes the proof of Theorem~\ref{thm:normconvergence}.

Perhaps the most interesting instance of Theorem~\ref{thm:normconvergence} is the convergence of \eqref{eq:pprodem} and \eqref{eq:pprodme} in the $\textup{L}^1$-norm assuming that $f$ and $g$ are arbitrary square-integrable functions, and one could even specialize $a=2$. Already for this particular case we do not see a simpler proof than repeating the above steps.

\section*{Acknowledgments}
This work is supported in part by the \emph{Croatian Science Foundation} under project UIP-2017-05-4129 (MUNHANAP). V. K. also acknowledges support of the \emph{Fulbright Scholar Program} and hospitality of the \emph{Georgia Institute of Technology} in the academic year 2019--2020. The authors are grateful to the anonymous reviewer for correcting and completing a couple of bibliographical references.



\begin{thebibliography}{99}

\bibitem{BB88:martpprod} R. Ba\~{n}uelos, A. G. Bennett, \emph{Paraproducts and commutators of martingale transforms}, Proc. Amer. Math. Soc. {\bf 103} (1988), no. 4, 1226--1234.

\bibitem{BMN10:pprod} \'{A}. B\'{e}nyi, D. Maldonado, V. Naibo, \emph{What is a Paraproduct?}, Notices Amer. Math. Soc. {\bf 57} (2010), no. 7, 858--860.

\bibitem{B31:pt} G. D. Birkhoff, \emph{Proof of the ergodic theorem}, Proc. Nat. Acad. Sci. U.S.A. {\bf 17} (1931), no. 12, 656--660.

\bibitem{Bur66:marttran} D. L. Burkholder, \emph{Martingale transforms}, Ann. Math. Statist. {\bf 37} (1966), 1494--1504.

\bibitem{Cal68:transf} A.-P. Calder\'{o}n, \emph{Ergodic theory and translation-invariant operators}, Proc. Nat. Acad. Sci. U.S.A. \textbf{59} (1968), 349--353.

\bibitem{CL92:martpprod} J.-A. Chao, R.-L. Long, \emph{Martingale transforms with unbounded multipliers}, Proc. Amer. Math. Soc. {\bf 114} (1992), no. 3, 831--838.

\bibitem{CL84:L2} J.-P. Conze, E. Lesigne, \emph{Th\'{e}or\`{e}mes ergodiques pour des mesures diagonales}, Bull. Soc. Math. France {\bf 112} (1984), no. 2, 143--175.

\bibitem{D40:mart} J. L. Doob, \emph{Regularity properties of certain families of chance variables}, Trans. Amer. Math. Soc. {\bf 47} (1940), 455--486.

\bibitem{DKST17:nvea} P. Durcik, V. Kova\v{c}, K. A. \v{S}kreb, C. Thiele, \emph{Norm variation of ergodic averages with respect to two commuting transformations}, Ergodic Theory Dynam. Systems {\bf 39} (2019), no. 3, 658--688.

\bibitem{D10:prob} R. Durrett, \emph{Probability: Theory and Examples}, fourth edition, Cambridge Series in Statistical and Probabilistic Mathematics, Cambridge University Press, 2010.

\bibitem{Fra16:problems} N. Frantzikinakis, \emph{Some open problems on multiple ergodic averages}, Bull. Hellenic Math. Soc. {\bf 60} (2016), 41--90.

\bibitem{FLW16:random} N. Frantzikinakis, E. Lesigne, M. Wierdl, \emph{Random sequences and pointwise convergence of multiple ergodic averages}, Indiana Univ. Math. J. {\bf 61} (2012), no. 2, 585--617.

\bibitem{FK78:msz} H. Furstenberg, Y. Katznelson, \emph{An ergodic Szemer\'{e}di theorem for commuting transformations}, J. Anal. Math. {\bf 38} (1978), no. 1, 275--291.

\bibitem{GLLZ:interp} L. Grafakos, L. Liu, S. Lu, F. Zhao, \emph{The multilinear Marcinkiewicz interpolation theorem revisited: the behavior of the constant}, J. Funct. Anal. {\bf 262} (2012), no. 5, 2289--2313.

\bibitem{G68:decomp} R. F. Gundy, \emph{A decomposition for {$L^{1}$}-bounded martingales}, Ann. Math. Statist. {\bf 39} (1968), 134--138.

\bibitem{IT63:open} A. Ionescu Tulcea, C. Ionescu Tulcea, \emph{Abstract ergodic theorems}, Trans. Amer. Math. Soc. {\bf 107} (1963), 107--124.

\bibitem{J59:open} M. Jerison, \emph{Martingale formulation of ergodic theorems}, Proc. Amer. Math. Soc. {\bf 10} (1959), 531--539.

\bibitem{JKRW98:Lp} R. L. Jones, R. Kaufman, J. M. Rosenblatt, M. Wierdl, \emph{Oscillation in ergodic theory}, Ergodic Theory Dynam. Systems {\bf 18} (1998), no. 4, 889--935.

\bibitem{JOR96:L2} R. L. Jones, I. V. Ostrovskii, J. M. Rosenblatt, \emph{Square functions in ergodic theory}, Ergodic Theory Dynam. Systems {\bf 16} (1996), no. 2, 267--305.

\bibitem{K98:open} A. G. Kachurovskii, \emph{A martingale ergodic theorem}, Mat. Zametki {\bf 64} (1998), no. 2, 311--314.

\bibitem{K07:open} A. G. Kachurovskii, \emph{General theories unifying ergodic averages and martingales}, Tr. Mat. Inst. Steklova {\bf 256} (2007), Din. Sist. i Optim., 172--200.

\bibitem{K52:open} S. Kakutani, \emph{Ergodic theory}, Proceedings of the International Congress of Mathematicians, Cambridge, 1950, vol. 2, pp. 128--142. Amer. Math. Soc., Providence, 1952.

\bibitem{K12:tp} V. Kova\v{c}, \emph{Boundedness of the twisted paraproduct}, Rev. Mat. Iberoam. {\bf 28} (2012), no. 4, 1143--1164.

\bibitem{KS17:mart} V. Kova\v{c}, K. A. \v{S}kreb, \emph{Bellman functions and $\textup{L}^p$ estimates for paraproducts}, Probab. Math. Statist. {\bf 38} (2018), no. 2, 459--479.

\bibitem{KS15:mart} V. Kova\v{c}, K. A. \v{S}kreb, \emph{One modification of the martingale transform and its applications to paraproducts and stochastic integrals}, J. Math. Anal. Appl. {\bf 426} (2015), no. 2, 1143--1163.

\bibitem{KZK19:mart} V. Kova\v{c}, P. Zorin-Kranich, \emph{Variational estimates for martingale paraproducts}, Electron. Commun. Probab. {\bf 24} (2019), paper no. 48, 14 pp.

\bibitem{N64:open} J. Neveu, \emph{Deux remarques sur la théorie des martingales}, Z. Wahrscheinlichkeitstheorie und Verw. Gebiete {\bf 3} (1964), 122--127. 

\bibitem{P84:open} D. Petz, \emph{Quantum ergodic theorems}, Quantum probability and applications to the quantum theory of irreversible processes (Villa Mondragone, 1982), 289--300, Lecture Notes in Math. {\bf 1055}, Springer, Berlin, 1984.

\bibitem{P10:open} I. V. Podvigin, \emph{A martingale-ergodic theorem}, Sibirsk. Mat. Zh. {\bf 51} (2010), no. 6, 1422--1429.

\bibitem{P14:open} I. V. Podvigin, \emph{Diagonal martingale ergodic sequences}, J. Math. Sci. (N.Y.) {\bf 198} (2014), no. 5, 602--607.

\bibitem{R61:open} G.-C. Rota, \emph{Une th\'{e}orie unifi\'{e}e des martingales et des moyennes ergodiques}, C. R. Acad. Sci. Paris {\bf 252} (1961), 2064--2066.

\bibitem{S20:open} F. A. Shahidi, Vector valued unified martingale and ergodic theorems with continuous parameter (2020), preprint, available at arXiv:2002.06399.

\bibitem{SG12:open} F. A. Shahidi, I. G. Ganiev, \emph{Vector valued martingale-ergodic and ergodic-martingale theorems}, Stoch. Anal. Appl. {\bf 30} (2012), no. 5, 916--932.


\bibitem{VK99:open} A. M. Vershik, A. G. Kachurovskii, \emph{Rates of convergence in ergodic theorems for locally finite groups, and reversed martingales}, Differ. Uravn. Protsessy Upr. (1999), no. 1, 19--26.

\bibitem{vN32:erg} J. von Neumann, \emph{Proof of the Quasi-Ergodic Hypothesis}, Proc. Nat. Acad. Sci. U.S.A. {\bf 18} (1932), no. 1, 70--82.

\end{thebibliography}
\end{document}